\documentclass[a4paper,11pt,leqno]{amsart}
\usepackage{amsmath,amsfonts,amssymb,amsthm}
\usepackage[a4paper, margin=2cm]{geometry}
\usepackage{hyperref}

\newtheorem{thm}{Theorem}[section]
\newtheorem{prop}[thm]{Proposition}
\newtheorem{coro}[thm]{Corollary}
\newtheorem{lem}[thm]{Lemma}
\theoremstyle{definition}
\newtheorem{defn}[thm]{Definition}
\newtheorem{rmk}[thm]{Remark}
\numberwithin{equation}{section}

\newcommand{\ord}{\operatorname{ord}}
\newcommand{\sgn}{\operatorname{sgn}}
\newcommand{\rank}{\operatorname{rank}}
\newcommand{\gen}{\operatorname{gen}}
\newcommand{\z}{\mathbb{Z}}

\newcommand{\Z}{\mathbb{Z}}
\newcommand{\R}{\mathbb{R}}
\newcommand{\mf}{\mathfrak}
\newcommand{\qf}[1]{\langle #1 \rangle}
\newcommand{\p}{\mf p}
\newcommand{\q}{\mf q}
\newcommand{\os}{\mf o}
\newcommand{\uu}{\mf u}
\newcommand{\nn}{\mf n}
\newcommand{\sca}{\mf s}
\newcommand{\vv}{\mf v}
\newcommand{\dd}{\mf d}
\newcommand{\DD}{\mf D}
\newcommand{\Hy}{\mathbb{H}}
\newcommand{\Ani}{\mathbb{A}}
\newcommand{\Psip}{\Psi_{\p}}

\newcommand{\refpart}[2]{\hyperref[#1]{\ref*{#1}#2}}

\begin{document}

\title[Regular quadratic lattices are rare over real quadratic fields]{Regular quadratic lattices are rare \\ over real quadratic fields}

\author{V\' \i t\v ezslav Kala}
\address{Charles University, Faculty of Mathematics and Physics, Department of Algebra, Sokolov\-sk\' a 83, 18600 Praha 8, Czech Republic}
\email{vitezslav.kala@matfyz.cuni.cz}

\author{Dayoon Park}
\address{Charles University, Faculty of Mathematics and Physics, Department of Algebra, Sokolov\-sk\' a 83, 18600 Praha 8, Czech Republic}
\email{pdy1016@snu.ac.kr}

\begin{abstract}
We show that a positive definite regular quadratic lattice of rank at least $4$ over a totally real number field contains a regular sublattice that represents precisely the totally positive integers lying in an
ideal. Under the generalized Riemann hypothesis and a conjecture on class numbers, it follows that, for each $m$, real quadratic fields admitting a positive definite regular lattice of some rank between $4$ and $m$
have density $0$. Unconditionally, only finitely many real quadratic fields of narrow class
number $1$ admit a positive definite regular lattice of some rank between $4$ and $7$.
\end{abstract}

\thanks{V.K. was supported by {Czech Science Foundation} grant 26-20514S. D.P. was supported by {Charles University} programme UNCE/24/SCI/022, and by the National Research Foundation of Korea (NRF) grant funded by the Korea government (MSIT) (No. RS-2025-02262988)}

\keywords{regular quadratic form, universal quadratic form, quadratic lattice, totally real number field, real quadratic field, Watson transformation}
\subjclass[2020]{11E08, 11E12, 11E20, 11R04, 11R11, 11R80}

\maketitle

\section{Introduction}\label{sec:intro}

A positive definite integral quadratic form over $\z$ is {\it regular} if it represents
every positive integer that is represented by its genus, i.e., every positive integer that it
represents over $\z_p$ for every prime $p$. Thus a regular form represents all positive integers
except those excluded by congruence conditions. Every {\it universal} form, i.e., every quadratic form
that represents all positive integers, is regular.

Regular forms were introduced by Dickson \cite{Di} in 1927, who also determined those of the type
$x^{2}+ay^{2}+bz^{2}$. Jones and Pall \cite{Jo, JP} then completed the
classification of regular diagonal ternary forms $ax^{2}+by^{2}+cz^{2}$. Watson \cite{Wa1}
proved that there are only finitely many equivalence classes of primitive positive definite
regular ternary forms by bounding their discriminants. Jagy--Kaplansky--Schiemann \cite{JKS}
produced a list of $913$ candidates for regular ternary forms and proved the regularity of
all but $22$ of them. Of these $22$ forms, $8$ were proved to be regular by Oh \cite{Oh2} and
$14$ by Lemke Oliver \cite{LO} under the generalized Riemann hypothesis. In rank at least $4$,
however, the analogous finiteness statement fails, since Earnest \cite{Ea1} constructed
infinitely many regular quaternary forms.

From now on, we will use the geometric language of quadratic lattices. Watson's proof of his finiteness theorem used transformations that replace a lattice by a
sublattice whose discriminant becomes smaller after rescaling 
\cite{Wa1,Wa2}. These transformations preserve
regularity and do not increase the class number, and by applying them repeatedly, one can reduce every regular ternary lattice
to one of bounded discriminant. Watson
transformations have since been used to determine regular ternary $\z$-lattices \cite{JKS,Oh2,LO},
to prove finiteness theorems also for $n$-regularity \cite{CO,Oh1} and
spinor regularity \cite{CE} (see \cite{CEO} for an overview), and to classify one-class genera
\cite{KL}.

Regular and universal quadratic lattices over the ring of integers $\os$ of a totally
real number field $F$ are defined analogously (see Section \ref{sec:notation}). Earnest
\cite{Ea2} considered regular lattices in this setting. Not too much is known about regular
ternary lattices over totally real number fields. Over each such field $F$, there are only
finitely many positive definite regular ternary $\os$-lattices up to isometry and rescaling.
Chan--Earnest--Icaza--Kim \cite{CEIK} proved this for $F=\mathbb Q(\sqrt5)$ by extending
Watson's method, and Chan--Icaza \cite{CI} for every $F$, using Watson transformations to study
ternary lattices that are {\it almost regular}, i.e., that represent all totally positive
integers represented by their genus, except for those in finitely many classes modulo squares of
units. All lattices in one-class genera are regular, and Pfeuffer \cite{Pf} showed that
positive definite one-class genera of rank at least $3$ exist over only finitely many totally real
number fields. Kirschmer--Lorch \cite{KL2} determined all the ternary ones, which form $4194$
genera over $30$ fields of degree at most $5$, and so at least these $30$ fields admit a positive
definite regular ternary lattice. It seems to be open whether only finitely many totally real
number fields do.

The main goal of this paper is to obtain the first results on the nonexistence of regular
lattices of large ranks, over a wide class of real quadratic fields. The principal idea is to
convert a regular lattice into one that is essentially universal (specifically, that represents
all the multiples of a fixed rational integer) and then to use the available results on the
nonexistence of such lattices. For that, the first step is to show that a regular lattice of
rank at least $4$ contains a regular sublattice that represents precisely the totally positive
elements of an ideal.

\begin{thm}\label{thm:main}
Let $F$ be a totally real number field with ring of integers $\os$ and let $L$ be a positive
definite regular $\os$-lattice of rank $n\ge4$. Then there is a regular $\os$-sublattice $L'$
of $L$ with $FL'=FL$ such that
\[
Q(L'\setminus\{0\})=2\sca(L')\cap\os^{+}.
\]
\end{thm}

Here $\os^{+}$ denotes the set of totally positive elements of $\os$. The scale $\sca(L')$ of $L'$
is the ideal of $\os$ generated by all values $B(v,w)$ of the bilinear form with $v,w\in L'$.

To prove Theorem \ref{thm:main}, we will apply Watson transformations successively at the
finitely many primes $\p$ of $F$ for which the localization $L_\p$ satisfies
$Q(L_\p)\ne2\sca(L_\p)$, until we obtain a lattice $L'$ with $Q(L'_\p)=2\sca(L'_\p)$ for every
prime $\p$. Over a general totally real field, the properties of these transformations that we
need (Lemma \ref{lem:E} and Propositions \ref{prop:index} and \ref{prop:shape}) are analogous
to the known ones over $\mathbb Q$, but most of them do not
seem to be available in the literature in this generality (see \cite{CI} for ternary lattices and
\cite{Be} for dyadic local fields).

To apply Theorem \ref{thm:main}, we will use universal quadratic lattices over totally real number fields and the known results that they are quite rare. Already
Siegel \cite{Si} showed that $\mathbb Q$ and $\mathbb Q(\sqrt5)$ are the only totally real
fields over which some sum of squares is universal, and 
Maa{\ss} \cite{Maa} proved that the sum of three squares is universal over $\mathbb Q(\sqrt5)$. By the asymptotic local--global principle of
Hsia--Kitaoka--Kneser \cite{HKK}, universal quadratic lattices exist over every totally real
number field, but their minimal rank is difficult to determine.

The main open problem is
Kitaoka's conjecture that only finitely many totally real fields admit a universal ternary
quadratic lattice. Over real quadratic fields, it was settled by
Chan--Kim--Raghavan \cite{CKR} for classical lattices and in general by Kim--Kim--Park
\cite{KKP}, who proved that only finitely many real quadratic fields admit a universal lattice
of rank at most $7$. Kala--Kr\'asensk\'y--Park--Yatsyna--\.Zmija \cite{KK+} then bounded the
discriminants of these fields explicitly. In arbitrary rank $r$, Blomer--Kala \cite{BK} first showed
that there are real quadratic fields admitting no universal quadratic lattice of rank at most
$r$, and Kala--Yatsyna--\.Zmija \cite{KYZ} even proved that real quadratic fields that do
admit such a lattice have density $0$. Over number fields of higher degrees, the minimal ranks of
universal quadratic lattices are known to be large in various families \cite{Kal,KT}. Man \cite{Man}
extended the density result of \cite{KYZ} to real multiquadratic fields.

In fact, \cite[Theorem 4.2]{KYZ} bounds the number of real quadratic fields that admit a lattice
of given rank representing every totally positive multiple of a fixed rational integer. By
Theorem \ref{thm:main}, every positive definite regular lattice of rank at least $4$ yields such
a lattice. Combining these with a conditional bound of
Bach--Sorenson \cite{BS} for the least prime ideal in a narrow ideal class, we obtain the following
result ($h^{+}(D)$ denotes the narrow class number of $\mathbb Q(\sqrt D)$).

\begin{thm}\label{thm:GRH}
Let $m\ge4$ be an integer. Assume the generalized Riemann hypothesis and assume that for every $\eta>0$ the set of squarefree integers $D>1$ with
$h^{+}(D)>D^{\eta}$ has density $0$ in the set of all squarefree integers. Then the set of
squarefree integers $D>1$ such that $\mathbb Q(\sqrt D)$ admits a positive definite regular
lattice of some rank $r$ with $4\le r\le m$ has density $0$ in the set of all squarefree
integers.
\end{thm}

The
assumption on $h^{+}(D)$ follows from a conjecture of Hooley \cite[Conjecture 7]{Hoo} on the average class
number of indefinite binary quadratic forms (see Remark \ref{rmk:hooley}). Under only a hypothesis
predicted by the Cohen--Lenstra heuristics, the same argument shows that there are $\gg X/\log X$
primes $p\le X$ for which $\mathbb Q(\sqrt p)$ admits no such lattice (Corollary \ref{cor:CL}).

If $L$ is a positive definite regular lattice of rank $n\ge4$ over a totally real number
field, then the orthogonal sum $L\perp\qf{a}$ is
regular for some $a\in\Z^+$ (Proposition \ref{prop:rankup}). Hence a
totally real number field that admits a positive definite regular lattice of rank $n\ge4$ admits
one of every rank $m\ge n$, so in Theorem \ref{thm:GRH} the condition ``of some rank $r$
with $4\le r\le m$'' is equivalent to ``of rank $m$''.

When the narrow class number of $F$ is $1$, a suitable rescaling of the lattice $L'$ from
Theorem \ref{thm:main} has scale $\os$ and represents exactly
the totally positive integers divisible by $2$ (Corollary \ref{cor:narrow}). In ranks at most $7$,
the explicit bounds of \cite{KK+} for such lattices give the following unconditional result.

\begin{thm}\label{thm:finite}
Let $F$ be a real quadratic field of narrow class number $1$ that admits a positive definite
regular lattice of rank $r$ for some $r$ with $4\le r\le7$. Then the discriminant of $F$ is
less than
\[
32\,400,\qquad 3\,240\,000,\qquad 2\,470\,090\,000,\qquad\text{and}\qquad
63\,234\,304\,000\,000
\]
for $r=4$, $5$, $6$, and $7$, respectively.
\end{thm}
Under the generalized
Riemann hypothesis, we also obtain $h^{+}(D)\gg D^{1/32}(\log D)^{-3/2}$ if $\mathbb Q(\sqrt D)$ admits a positive definite
regular lattice of some rank between $4$ and $7$ (Corollary \ref{cor:GRHclass}).

Theorem \ref{thm:main} holds over
every totally real number field, so a counting result of the type of \cite{KYZ} in higher degrees,
together with the analogous assumption on narrow class numbers, would extend
Theorem \ref{thm:GRH} to those degrees. However, the proof in \cite{KYZ} uses continued fractions
of quadratic irrationalities, which have no analogue in all fields of a given degree that would be useful for this purpose.
An extension to real multiquadratic fields using \cite{Man} may be possible, though.

\section*{Acknowledgments}

We thank the members of the \href{https://www.karlin.mff.cuni.cz/~kala/web/ufoclan}{UFOCLAN} research group for our helpful discussions.

The LLMs Claude Opus 5, 5.5, and Fable 5.2 prepared the text of the manuscript from an initial draft written by the authors.
The ideas behind the results presented in the paper and their proofs are primarily human-generated. A few notable exceptions suggested by AI are the present uniform proof of Proposition \ref{prop:red}, which streamlines the authors' original argument, and some of the specific conjectures used in Section \ref{sec:density}. The authors maintain full responsibility for the contents of the paper and their correctness.

\section{Preliminaries}\label{sec:notation}

\subsection*{Number fields}

Throughout, $F$ is a totally real number field of degree $d$ over $\mathbb Q$ with
ring of integers $\os$. We denote by $\sigma_1,\dots,\sigma_d\colon F\to\R$ the embeddings of $F$
into $\R$. An element $\alpha\in F$ is {\it totally positive} if $\sigma_i(\alpha)>0$ for all
$1\le i\le d$. We write $F^{+}$ for the set of totally positive elements of $F$ and
$\os^{+}:=\os\cap F^{+}$ for the set of totally positive integers of $F$. By a {\it prime} of
$F$ we mean a nonzero prime ideal of $\os$.

Two fractional ideals $\mf a$ and $\mf b$ of $F$ are {\it narrowly equivalent} if
$\mf b=\alpha\mf a$ for some totally positive $\alpha\in F$. The equivalence classes are called
{\it narrow ideal classes} and they form a finite abelian group, the
{\it narrow ideal class group} of $F$, whose order is the {\it narrow class number} of $F$.

Let $\p$ be a prime of $F$. We denote by $F_\p$ the completion of $F$ at $\p$, by $\os_\p$ its
ring of integers, and by $\uu_\p$ the group of units of $\os_\p$. We also write $\p$ for the
maximal ideal of $\os_\p$, fix a generator $\pi_\p$ of it, and let
$\ord_\p\colon F_\p\to\z\cup\{\infty\}$ be the valuation with $\ord_\p(\pi_\p)=1$. For a
nonzero fractional ideal $\mf a$ of $F$ or of $F_\p$, we write $\ord_\p\mf a$ for the integer
$k$ with $\mf a\os_\p=\p^{k}$. The prime $\p$ is {\it dyadic} if $2\in\p$. For subsets $X,Y$ of
$F$ or of $F_\p$ and an element $\alpha$ of $F$ or of $F_\p$, we write
\[
X+Y:=\{x+y\mid x\in X,\ y\in Y\}\qquad\text{and}\qquad \alpha X:=\{\alpha x\mid x\in X\}.
\]
We also put $\os_\p^{2}:=\{x^{2}\mid x\in\os_\p\}$,
$F_\p^{\times2}:=\{x^{2}\mid x\in F_\p^{\times}\}$, and
$\uu_\p\os_\p^{2}:=\{\varepsilon x^{2}\mid \varepsilon\in\uu_\p,\ x\in\os_\p\}$.

For $\xi\in F_\p$, the {\it quadratic defect} of $\xi$ is
$\DD(\xi):=\bigcap_{\eta\in F_\p}(\xi-\eta^{2})\os_\p$ (see \cite[63A]{O}). For every dyadic
prime $\p$ we fix a unit $\rho_\p\in\uu_\p$ such that
\[
\Delta_\p:=1-4\rho_\p
\]
has quadratic defect $\DD(\Delta_\p)=4\os_\p$. Such a unit exists and $F_\p(\sqrt{\Delta_\p})$
is the unramified quadratic extension of $F_\p$ (see \cite[63:2--63:4]{O}).

\subsection*{Quadratic spaces and lattices}

Let $E$ be one of the fields $F$, $F_\p$, and $\R$, where $\p$ is a prime of $F$, and let $V$ be
a finite-dimensional vector space over $E$. If $B$ is a nondegenerate symmetric bilinear form
on $V$, then we call $V$ (together with $B$) a {\it quadratic space} over $E$ and $Q(v):=B(v,v)$
its {\it quadratic form}. For a subset $X\subset V$ we write $Q(X):=\{Q(v)\mid v\in X\}$. A
nonzero vector $v\in V$ is {\it isotropic} if $Q(v)=0$. The space $V$ is {\it isotropic} if it
contains an isotropic vector and {\it anisotropic} otherwise. The {\it Gram matrix} of a tuple $v_1,\dots,v_n\in V$ is the matrix $\big(B(v_i,v_j)\big)_{1\le i,j\le n}$. The
{\it orthogonal sum} $V_1\perp V_2$ of two quadratic spaces $V_1$ and $V_2$ over $E$, with
bilinear forms $B_1$ and $B_2$, is the vector space $V_1\oplus V_2$ with the bilinear form
$B\big((v_1,v_2),(w_1,w_2)\big):=B_1(v_1,w_1)+B_2(v_2,w_2)$.

Now let $E$ be $F$ or $F_\p$, let $\os_E$ be its ring of integers, and let $V$ be a quadratic space over $E$. A {\it lattice on $V$} is a finitely
generated $\os_E$-submodule $L$ of $V$ with $EL=V$. We then call $L$ an {\it $\os_E$-lattice}, or a {\it lattice over} $E$,
and say that $V$ is the quadratic space {\it spanned} by $L$. The {\it rank} of $L$ is
$\rank L:=\dim V$. We write $B(L,L):=\{B(v,w)\mid v,w\in L\}$. Unless stated otherwise, all
lattices are assumed to be {\it integral}, i.e., $B(L,L)\subset\os_E$. Integral lattices in this
sense are also called {\it classical}, e.g., in \cite{KYZ}. A lattice is {\it isotropic} if the
quadratic space spanned by it is isotropic, and {\it anisotropic} otherwise.

Two lattices $L$ and $L'$ on quadratic spaces $V$ and $V'$ with bilinear forms $B$ and $B'$ are
{\it isometric}, $L\cong L'$, if there is a linear bijection $\varphi\colon V\to V'$ with
$\varphi(L)=L'$ and $B'(\varphi(v),\varphi(w))=B(v,w)$ for all $v,w\in V$. The
{\it orthogonal sum} $L_1\perp L_2$ of lattices $L_1$ and $L_2$ on $V_1$ and $V_2$ is the lattice
$L_1\oplus L_2$ on $V_1\perp V_2$. If $L_1$ and $L_2$ are submodules of a lattice $L$ with
$L=L_1\oplus L_2$ and $B(L_1,L_2)=0$, then we also write $L=L_1\perp L_2$ and call this an
{\it orthogonal splitting} of $L$.

A {\it basis} of an $\os_E$-lattice $L$ is a basis of $L$ as an $\os_E$-module, which exists
only if $L$ is free.
For nonzero $a_1,\dots,a_n\in\os_E$, we write $\qf{a_1,\dots,a_n}$ for the lattice with a basis
$v_1,\dots,v_n$ such that $Q(v_i)=a_i$ for all $i$ and $B(v_i,v_j)=0$ for all $i\ne j$. For
$\alpha,\beta\in\os_\p$ with $\alpha\beta\ne1$, we write $A(\alpha,\beta)$ for the
$\os_\p$-lattice with a basis $v_1,v_2$ such that $Q(v_1)=\alpha$, $B(v_1,v_2)=1$, and
$Q(v_2)=\beta$, i.e., with Gram matrix
\[
\begin{pmatrix}\alpha&1\\1&\beta\end{pmatrix}.
\]

The {\it norm} $\nn(L)$ and the {\it scale} $\sca(L)$ of a lattice $L$ are the fractional ideals
of $E$ generated by $Q(L)$ and by $B(L,L)$, respectively. Since
$2B(v,w)=Q(v+w)-Q(v)-Q(w)$ and $Q(v)=B(v,v)$, we have
\begin{equation}\label{eq:nsn}
2\sca(L)\subset\nn(L)\subset\sca(L).
\end{equation}
Every lattice $L$ can be written as $L=\mf a_1v_1+\dots+\mf a_nv_n$ with a basis
$v_1,\dots,v_n$ of $V$ and fractional ideals $\mf a_1,\dots,\mf a_n$ of $E$ (see
\cite[\S81]{O}). The {\it volume} of $L$ is the fractional ideal
\[
\vv(L):=\mf a_1^{2}\cdots\mf a_n^{2}\det\big(B(v_i,v_j)\big),
\]
which does not depend on these choices (see \cite[\S82]{O}). If $E=F_\p$, then $\os_\p$ is a principal ideal domain, and so $L$ has a basis
$v_1,\dots,v_n$ and
$\vv(L)=\det\big(B(v_i,v_j)\big)\os_\p$. For lattices $L$ and $L'$ we have
\begin{equation}\label{eq:perp}
\begin{gathered}
Q(L\perp L')=Q(L)+Q(L'),\quad \nn(L\perp L')=\nn(L)+\nn(L'),\\
\sca(L\perp L')=\sca(L)+\sca(L'),\quad \vv(L\perp L')=\vv(L)\vv(L').
\end{gathered}
\end{equation}

Let $K\subset L$ be lattices on $V$. By the structure theorem for finitely generated torsion
modules over Dedekind domains, $L/K\cong\os_E/\mf c_1\oplus\dots\oplus\os_E/\mf c_k$ for some
nonzero ideals $\mf c_1,\dots,\mf c_k$ of $\os_E$. The {\it index ideal} of $K$ in $L$ is
$[L:K]:=\mf c_1\cdots\mf c_k$, which does not depend on this decomposition. By \cite[82:11]{O}
we have $\vv(K)=[L:K]^{2}\vv(L)$. A vector $v$ of an $\os_\p$-lattice $L$ is {\it primitive} if
$v\notin\p L$.

For $\alpha\in E^{\times}$, we denote by $L^{\alpha}$ the lattice $L$ on the quadratic space $V$
with the bilinear form $\alpha B$ in place of $B$. This should not be confused with the module
$\alpha L$, which is a lattice on $V$ with the bilinear form $B$. With $n=\rank L$, we have
\begin{equation}\label{eq:scaling}
\begin{gathered}
Q(L^{\alpha})=\alpha Q(L),\quad \nn(L^{\alpha})=\alpha\nn(L),\quad
\sca(L^{\alpha})=\alpha\sca(L),\quad \vv(L^{\alpha})=\alpha^{n}\vv(L),\\
Q(\alpha L)=\alpha^{2}Q(L),\quad \nn(\alpha L)=\alpha^{2}\nn(L),\quad
\sca(\alpha L)=\alpha^{2}\sca(L),\quad \vv(\alpha L)=\alpha^{2n}\vv(L).
\end{gathered}
\end{equation}

An $\os_\p$-lattice $L$ of rank $n$ is called {\it $\sca(L)$-modular}, or simply {\it modular},
if $\vv(L)=\sca(L)^{n}$. An $\os_\p$-modular lattice is also called {\it unimodular}. By
\eqref{eq:scaling},
\[
L \text{ is } \p^{i}\text{-modular}\iff L^{\pi_\p^{-i}} \text{ is unimodular}.
\]
Every $\os_\p$-lattice $L$ has a {\it Jordan splitting}, i.e., an orthogonal splitting
\[
L=L_1\perp L_2\perp\dots\perp L_j
\]
in which each $L_i$ is modular and $\sca(L_1)\supsetneq\sca(L_2)\supsetneq\dots\supsetneq\sca(L_j)$
(see \cite[91C]{O}), we call $L_1$
a {\it leading component} of $L$. A Jordan splitting is in general not unique, so every
statement below about a leading component refers to a fixed Jordan splitting.

\subsection*{Localization and representations}

Let $L$ be an $\os$-lattice on a quadratic space $V$ over $F$. For a prime $\p$ of $F$, we view
$V$ as a subset of the quadratic space $V_\p:=V\otimes_F F_\p$ over $F_\p$. The
{\it localization} of $L$ at $\p$ is the $\os_\p$-lattice $L_\p:=\os_\p L$ on $V_\p$. Then $L$ is
dense in $L_\p$ and $\nn(L)\os_\p=\nn(L_\p)$, $\sca(L)\os_\p=\sca(L_\p)$, and
$\vv(L)\os_\p=\vv(L_\p)$ (see \cite[81C]{O}). For each embedding $\sigma_i$, we denote by
$V_{\sigma_i}:=V\otimes_{F,\sigma_i}\R$ the real quadratic space obtained from $V$ via $\sigma_i$
and by $Q_{\sigma_i}$ its quadratic form. The lattice $L$ is {\it positive definite} if
$Q_{\sigma_i}(v)>0$ for all $1\le i\le d$ and all nonzero $v\in V_{\sigma_i}$.

An element $\alpha\in\os$ is {\it represented} by $L$ if $\alpha\in Q(L)$. It is
{\it locally represented} by $L$ if $\alpha\in Q(L_\p)$ for every prime $\p$ of $F$ and
$\sigma_i(\alpha)\in Q_{\sigma_i}(V_{\sigma_i})$ for all $1\le i\le d$. We denote by $Q(\gen L)$ the
set of elements of $\os$ that are locally represented by $L$. The notation refers to the
{\it genus} $\gen L$ of $L$, i.e., the set of lattices $K$ on $V$ with $K_\p\cong L_\p$ for every
prime $\p$, since $Q(\gen L)$ is the union of the sets $Q(K)$ with $K\in\gen L$ (see
\cite[102:5]{O}).
Clearly $Q(L)\subset Q(\gen L)$. We call
\[
E(L):=Q(\gen L)\setminus Q(L)
\]
the {\it exceptional set} of $L$.
The lattice $L$ is {\it regular} if $E(L)=\emptyset$, i.e., if $Q(L)=Q(\gen L)$. A positive
definite lattice $L$ is {\it universal} if $\os^{+}\subset Q(L)$.

\begin{lem}\label{lem:units}
Let $\p$ be a prime of $F$, let $\alpha\in F_\p^{\times}$, and let $\varepsilon\in\uu_\p$. Then:
\begin{itemize}
\item[(1)] $\alpha\uu_\p\os_\p^{2}
=\{x\in F_\p\mid \ord_\p(x)\in\ord_\p(\alpha)+2\z_{\ge0}\}\cup\{0\}$.
\item[(2)] $\alpha\uu_\p\os_\p^{2}+\alpha\pi_\p\uu_\p\os_\p^{2}=\alpha\os_\p$.
\item[(3)] $\alpha\uu_\p\os_\p^{2}+\varepsilon\alpha\os_\p^{2}=\alpha\os_\p$.
\end{itemize}
\end{lem}

\begin{proof}
(1) If $x=\alpha\eta y^{2}$ with $\eta\in\uu_\p$ and $y\in\os_\p\setminus\{0\}$, then
$\ord_\p(x)=\ord_\p(\alpha)+2\ord_\p(y)$. Conversely, if $\ord_\p(x)=\ord_\p(\alpha)+2c$ with
$c\ge0$, then $x=\alpha\eta(\pi_\p^{c})^{2}$ with the unit $\eta:=x\alpha^{-1}\pi_\p^{-2c}$.

(2) Both sets on the left are contained in $\alpha\os_\p$ and contain $0$, so their sum
contains their union, which is $\alpha\os_\p$ by (1).

(3) Both sets on the left are contained in $\alpha\os_\p$, hence so is their sum. Conversely,
let $z\in\alpha\os_\p$ be nonzero. If $\ord_\p(z)-\ord_\p(\alpha)$ is even, then
$z\in\alpha\uu_\p\os_\p^{2}$ by (1), so $z=z+\varepsilon\alpha\cdot0^{2}$ lies in the sum. If
$\ord_\p(z)-\ord_\p(\alpha)$ is odd, then $\ord_\p(z)>\ord_\p(\alpha)$, hence
$\ord_\p(z-\varepsilon\alpha)=\ord_\p(\alpha)$. Thus $z-\varepsilon\alpha\in\alpha\uu_\p
\os_\p^{2}$ by (1), so $z=(z-\varepsilon\alpha)+\varepsilon\alpha\cdot1^{2}$ lies in the sum.
\end{proof}

\begin{lem}\label{lem:d}
Let $L$ be a lattice of rank $n$ over $\os$ or over $\os_\p$. Then
$\vv(L)\subset\sca(L)^{n}$ and the ideal
\[
\dd(L):=2^{n}\nn(L)^{-n}\vv(L)
\]
is integral. If $L$ is an $\os$-lattice, then moreover $\dd(L)\os_\p=\dd(L_\p)$ for every prime
$\p$ of $F$.
\end{lem}

\begin{proof}
The inclusion $\vv(L)\subset\sca(L)^{n}$ is \cite[82:10]{O}. By \eqref{eq:nsn} we have
$2\sca(L)\subset\nn(L)$, so $\vv(L)\subset\sca(L)^{n}\subset 2^{-n}\nn(L)^{n}$, i.e., $\dd(L)$
is integral. The last assertion follows from
$\nn(L)\os_\p=\nn(L_\p)$ and $\vv(L)\os_\p=\vv(L_\p)$.
\end{proof}

\begin{lem}\label{lem:scaleE}
Let $L$ be a lattice on a quadratic space over $F$ and let $\alpha\in F^{\times}$ be such that
$L^{\alpha}$ is integral. Then $Q(\gen L^{\alpha})=\alpha Q(\gen L)$ and
$E(L^{\alpha})=\alpha E(L)$. In particular, if $L$ is regular then so is $L^{\alpha}$.
\end{lem}

\begin{proof}
By \eqref{eq:scaling}, $Q(L^{\alpha}_\p)=\alpha Q(L_\p)\subset\os_\p$ for every prime $\p$ of $F$
and $Q_{\sigma_i}(V_{\sigma_i}^{\alpha})=\sigma_i(\alpha)Q_{\sigma_i}(V_{\sigma_i})$ for every
$i$. Hence $\alpha Q(\gen L)\subset\os$ and an element $\beta\in\os$ is locally represented
by $L^{\alpha}$ if and only if $\alpha^{-1}\beta$ is locally represented by $L$, which gives
the first identity. The second follows from $Q(L^{\alpha})=\alpha Q(L)$.
\end{proof}

\section{Modular lattices}\label{sec:modular}

In this section we will describe, for a modular $\os_\p$-lattice $L$, the set $Q(L)$ and the
ideals $Q(v)\os_\p$ for the primitive vectors $v\in L$. Most results below are classical (see \cite[\S\S92--93]{O}), but we include them for completeness.

For anisotropic binary lattices,
at non-dyadic and at dyadic primes alike, we will use the following lemma.

\begin{lem}\label{lem:unram}
Let $\p$ be a prime of $F$, let $E$ be the unramified quadratic extension of $F_\p$, and let
$\theta$ be an element of its ring of integers $\os_E$ such that $\os_E=\os_\p[\theta]$. Then:
\begin{itemize}
\item[(1)] $N_{E/F_\p}(x+y\theta)\in\uu_\p$ for all $x,y\in\os_\p$ that do not both lie in $\p$.
\item[(2)] $\{N_{E/F_\p}(x+y\theta)\mid x,y\in\os_\p\}=\uu_\p\os_\p^{2}$.
\end{itemize}
\end{lem}

\begin{proof}
As $E/F_\p$ is unramified, $\p\os_E$ is the maximal ideal of $\os_E$ and
$N_{E/F_\p}(E^{\times})=\{x\in F_\p^{\times}\mid \ord_\p(x)\in2\z\}$ (see
\cite[Ch. V, \S2, Corollary to Proposition 3 and Remark 1]{Se}).

(1) Since $1,\theta$ is a basis of $\os_E$ over $\os_\p$, the element $x+y\theta$ does not lie
in $\p\os_E$. Hence it is a unit of $\os_E$ and its norm is a unit of $\os_\p$.

(2) The set on the left is $N_{E/F_\p}(\os_E)$. A nonzero element of $E$ lies in $\os_E$ if
and only if its norm lies in $\os_\p$, so
$N_{E/F_\p}(\os_E)=\{x\in F_\p^{\times}\mid \ord_\p(x)\in2\z_{\ge0}\}\cup\{0\}$, which is
$\uu_\p\os_\p^{2}$ by Lemma \refpart{lem:units}{(1)}.
\end{proof}

\subsection*{Non-dyadic primes}

\begin{prop}\label{prop:nondy}
Let $\p$ be a non-dyadic prime and let $L$ be a nonzero $\p^{i}$-modular $\os_\p$-lattice of
rank $n$. Then:
\begin{itemize}
\item[(1)] If $n\ge3$, then $L$ is isotropic.
\item[(2)] If $L$ is isotropic, then $Q(L)=\p^{i}$.
\item[(3)] If $L$ is anisotropic, then $n\le 2$, $Q(v)\os_\p=\p^{i}$ for every primitive $v\in
L$, and $Q(L)\subset\pi_\p^{i}\uu_\p\os_\p^{2}$, with equality if $n=2$.
\end{itemize}
\end{prop}

\begin{proof}
	After replacing $L$ by the unimodular lattice $L^{\pi_\p^{-i}}$, we may assume that $i=0$.
	
	(1) Let $n\ge3$. As $L$ has an orthogonal basis (see \cite[\S92]{O}), we can write
	$L=\os_\p w\perp L'$ with $Q(w)\in\uu_\p$ and a unimodular lattice $L'$ of rank $n-1\ge2$. By
	\cite[92:1b]{O}, $L'$ represents every element of $\uu_\p$, in particular $-Q(w)$, so $L$ is
	isotropic.
	
	(2) Since $L$ is isotropic, $n\ge2$. If $n\ge3$, then $Q(L)=\os_\p$ by \cite[92:1b]{O}. If
	$n=2$, then $L\cong\qf{1,\varepsilon}$ for some $\varepsilon\in\uu_\p$ by \cite[92:1]{O}. As
	$L$ is isotropic, $-\varepsilon$ is a square in $\uu_\p$, so $L\cong\qf{1,-1}$. This lattice
	represents every $\alpha\in\os_\p$, since
	$\alpha=\big(\frac{\alpha+1}2\big)^{2}-\big(\frac{\alpha-1}2\big)^{2}$ with
	$\frac{\alpha\pm1}2\in\os_\p$. Hence $Q(L)=\os_\p$.
	
	(3) That $n\le2$ follows from (1). If $n=1$, then $L=\os_\p w$ with $Q(w)\in\uu_\p$, the
	primitive vectors of $L$ are the vectors $xw$ with $x\in\uu_\p$, and $Q(L)=Q(w)\os_\p^{2}\subset
	\uu_\p\os_\p^{2}$. Let $n=2$. By \cite[92:1]{O}, $L=\os_\p w_1\perp\os_\p w_2$ with $Q(w_1)=1$
	and $Q(w_2)=\varepsilon$, where $\varepsilon\in\uu_\p$ is the determinant of a Gram matrix of
	$L$. As $L$ is anisotropic, $-\varepsilon$ is not a square in $F_\p$, so
	$E:=F_\p(\sqrt{-\varepsilon})$ is the unramified quadratic extension of $F_\p$. Its ring of
	integers is $\os_E=\os_\p[\sqrt{-\varepsilon}]$ and
	\[
	Q(xw_1+yw_2)=x^{2}+\varepsilon y^{2}=N_{E/F_\p}(x+y\sqrt{-\varepsilon})\qquad(x,y\in\os_\p).
	\]
	The vector $xw_1+yw_2$ is primitive if and only if $x$ and $y$ do not both lie in $\p$. Hence
	Lemma \ref{lem:unram} gives $Q(v)\in\uu_\p$ for every primitive $v\in L$ and
	$Q(L)=\uu_\p\os_\p^{2}$.
\end{proof}

\begin{coro}\label{cor:nondy4}
Let $\p$ be a non-dyadic prime, let $M$ be a $\p^{a}$-modular and $N$ a
$\p^{a+1}$-modular $\os_\p$-lattice, both of rank $2$. Then $Q(M\perp N)=\p^{a}$.
\end{coro}

\begin{proof}
By Proposition \refpart{prop:nondy}{(2),(3)} we have $\pi_\p^{a}\uu_\p\os_\p^{2}\subset Q(M)$ and
$\pi_\p^{a+1}\uu_\p\os_\p^{2}\subset Q(N)$, so $\p^{a}\subset Q(M)+Q(N)=Q(M\perp N)$ by
Lemma \refpart{lem:units}{(2)} and \eqref{eq:perp}. The reverse inclusion follows from $Q(M\perp N)\subset
\nn(M\perp N)=\p^{a}$.
\end{proof}

\subsection*{Dyadic primes}

In the dyadic case we will only need modular lattices $L$ with $\nn(L)=2\sca(L)$, because
lattices with $\nn(L)\supsetneq2\sca(L)$ will be handled directly in
Proposition \ref{prop:index}.

\begin{prop}\label{prop:binary}
Let $\p$ be a dyadic prime. Up to isometry there are exactly two binary unimodular
$\os_\p$-lattices of norm $2\os_\p$, namely
\[
\Hy:=A(0,0)\qquad\text{and}\qquad \Ani:=A(2,2\rho_\p).
\]
The lattice $\Hy$ is isotropic and $Q(\Hy)=2\os_\p$. The lattice $\Ani$ is anisotropic,
$Q(\Ani)=2\uu_\p\os_\p^{2}$, and $Q(v)\os_\p=2\os_\p$ for every primitive $v\in\Ani$.
\end{prop}

\begin{proof}
The classification is \cite[93:11]{O}. Let $v_1,v_2$ be a basis of $\Hy$ with
$Q(v_1)=Q(v_2)=0$ and $B(v_1,v_2)=1$. Then $Q(xv_1+yv_2)=2xy$ for all $x,y\in\os_\p$, so $\Hy$ is
isotropic and $Q(\Hy)=2\os_\p$.

Let $v_1,v_2$ be a basis of $\Ani$ with $Q(v_1)=2$, $B(v_1,v_2)=1$, and $Q(v_2)=2\rho_\p$. Then
\[
Q(xv_1+yv_2)=2(x^{2}+xy+\rho_\p y^{2})\qquad(x,y\in\os_\p).
\]
Let $E:=F_\p(\sqrt{\Delta_\p})$, the unramified quadratic extension of $F_\p$, and let
$\omega:=(1+\sqrt{\Delta_\p})/2$, a root of $T^{2}-T+\rho_\p$. The discriminant of this
polynomial is the unit $\Delta_\p$, so $\os_E=\os_\p[\omega]$. We have
$N_{E/F_\p}(x+y\omega)=x^{2}+xy+\rho_\p y^{2}$ for all $x,y\in F_\p$, and the norm of a nonzero
element of $E$ is nonzero, so $\Ani$ is anisotropic. The vector $xv_1+yv_2$ is primitive if
and only if $x$ and $y$ do not both lie in $\p$. Hence Lemma \ref{lem:unram} gives
$Q(\Ani)=2\uu_\p\os_\p^{2}$ and $Q(v)\in2\uu_\p$ for every primitive $v\in\Ani$.
\end{proof}

\begin{prop}\label{prop:dymodular}
Let $\p$ be a dyadic prime and let $L$ be a nonzero $\p^{i}$-modular $\os_\p$-lattice with
$\nn(L)=2\sca(L)$. Then $\rank L=2r$ is even and
\[
L\cong \underbrace{\Hy^{\pi_\p^{i}}\perp\dots\perp\Hy^{\pi_\p^{i}}}_{r}
\qquad\text{or}\qquad
L\cong \underbrace{\Hy^{\pi_\p^{i}}\perp\dots\perp\Hy^{\pi_\p^{i}}}_{r-1}\perp\Ani^{\pi_\p^{i}}.
\]
Moreover, the following hold:
\begin{itemize}
\item[(1)] $2\pi_\p^{i}\uu_\p\os_\p^{2}\subset Q(L)$.
\item[(2)] $L$ is anisotropic if and only if $L\cong\Ani^{\pi_\p^{i}}$, in which case
$Q(L)=2\pi_\p^{i}\uu_\p\os_\p^{2}$ and $Q(v)\os_\p=\nn(L)$ for every primitive $v\in L$.
\item[(3)] If $L$ is isotropic, then $Q(L)=2\p^{i}=\nn(L)$.
\end{itemize}
\end{prop}

\begin{proof}
After replacing $L$ by the unimodular lattice $L^{\pi_\p^{-i}}$, we may assume that $i=0$, so
that $\nn(L)=2\os_\p\subsetneq\os_\p=\sca(L)$. By \cite[93:15]{O}, $L$ is an orthogonal sum of binary sublattices
$K_1,\dots,K_r$, which are unimodular since $L$ is, and so $\rank L=2r$. For each $s$ we have
$2\os_\p=2\sca(K_s)\subset\nn(K_s)\subset\nn(L)=2\os_\p$ by \eqref{eq:nsn}, so
$\nn(K_s)=2\os_\p$ and $K_s\cong\Hy$ or $K_s\cong\Ani$ by Proposition \ref{prop:binary}. Since
$\Ani\perp\Ani\cong\Hy\perp\Hy$ by \cite[93:18(ii)]{O}, we may replace pairs of summands
isometric to $\Ani$ by pairs isometric to $\Hy$, which gives the two stated normal forms.

(1) By Proposition \ref{prop:binary}, $Q(K_s)\supset2\uu_\p\os_\p^{2}$ for every $s$.

(2) If $L\cong\Ani$, then $L$ is anisotropic and the assertions are those of
Proposition \ref{prop:binary}. Otherwise $L$ contains a sublattice isometric to $\Hy$ and is
therefore isotropic.

(3) If $L$ is isotropic, then $L\not\cong\Ani$ by (2), so at least one summand in the normal
form above is isometric to $\Hy$ and
$2\os_\p=Q(\Hy)\subset Q(L)\subset\nn(L)=2\os_\p$.
\end{proof}

\begin{coro}\label{cor:dycoro}
Let $\p$ be a dyadic prime and let $L$ be an $\os_\p$-lattice with
$\nn(L)=2\sca(L)=2\p^{i}$. Let $L=M\perp N$ where $M$ is $\p^{i}$-modular with $\nn(M)=\nn(L)$.
If at least one of
\begin{itemize}
\item[(1)] $M$ is isotropic,
\item[(2)] $\nn(N)=\nn(M)$,
\item[(3)] $Q(N)\supset2\pi_\p^{i+1}\uu_\p\os_\p^{2}$
\end{itemize}
holds, then $Q(L)=\nn(L)=2\sca(L)$.
\end{coro}

\begin{proof}
Since $Q(L)\subset\nn(L)=2\p^{i}$ and $Q(L)=Q(M)+Q(N)$ by \eqref{eq:perp}, it suffices in each
case to show that $Q(M)+Q(N)\supset2\p^{i}$. Proposition \ref{prop:dymodular} applies to $M$,
since $\nn(M)=2\p^{i}=2\sca(M)$.

(1) Here $Q(M)=2\p^{i}$ by Proposition \refpart{prop:dymodular}{(3)}.

(2) Since $\nn(N)=2\p^{i}$, there is $w\in N$ with $Q(w)=2\varepsilon\pi_\p^{i}$ for some
$\varepsilon\in\uu_\p$, so $Q(N)\supset 2\varepsilon\pi_\p^{i}\os_\p^{2}$. As
$Q(M)\supset2\pi_\p^{i}\uu_\p\os_\p^{2}$ by Proposition \refpart{prop:dymodular}{(1)},
Lemma \refpart{lem:units}{(3)} gives $Q(M)+Q(N)\supset2\p^{i}$.

(3) Here $Q(M)\supset2\pi_\p^{i}\uu_\p\os_\p^{2}$ by Proposition \refpart{prop:dymodular}{(1)}, so
Lemma \refpart{lem:units}{(2)} gives $Q(M)+Q(N)\supset2\p^{i}$.
\end{proof}

\section{The transformation $\Psip$}\label{sec:psi}

In this section we will recall the transformation $\Psip$ defined by Chan--Icaza \cite{CI},
which will play the role of Watson's transformations in our proofs, and establish its basic
properties. Ternary analogues of many statements below are also due to them \cite[Section 3]{CI}.

\begin{defn}\label{def:psi}
Let $\p$ be a prime of $F$ and let $L$ be an $\os$-lattice or an $\os_\p$-lattice. We set
\[
\Psip(L):=\{v\in L\mid Q(v)\in\p\cdot\nn(L)\}.
\]
We put $\Psip^{0}(L):=L$ and, whenever $t\ge1$ and $\Psip^{t-1}(L)$ is a lattice,
$\Psip^{t}(L):=\Psip\big(\Psip^{t-1}(L)\big)$.
\end{defn}

The set $\Psip(L)$ need not be a lattice. For example, if $p$ is an odd rational prime and
$L=\z v_1\perp\z v_2\perp\z v_3=\qf{1,1,p-1}$, then $v_1+v_3$ and $v_2+v_3$ lie in
$\Psi_{p}(L)$, but their difference $v_1-v_2$ does not.

When $F=\mathbb Q$, the transformation $\Psip$ is closely related to Watson's transformations. For
a $\z$-lattice $L$ with $\nn(L)=2\z$ and an odd prime $p$, the sublattice
$\Lambda_{2p}(L)=\{x\in L\mid Q(x+z)\equiv Q(z)\pmod {2p}\text{ for all }z\in L\}$ of
\cite{Wa2} is contained in $\Psi_p(L)$, with equality if a leading component of
$L_p$ is anisotropic (see Propositions \refpart{prop:nondy}{(3)} and \ref{prop:shape}).

\begin{lem}\label{lem:psibasic}
Let $L$ be an $\os$-lattice on $V$ and let $\p$ be a prime of $F$. Then:
\begin{itemize}
\item[(1)] $\p L\subset\Psip(L)\subsetneq L$.
\item[(2)] $\Psip(L)=L\cap\Psip(L_\p)$.
\item[(3)] If $\Psip(L_\p)$ is an $\os_\p$-lattice, then $\Psip(L)$ is an $\os$-lattice on $V$
with $\Psip(L)_\p=\Psip(L_\p)$ and $\Psip(L)_\q=L_\q$ for every prime $\q\ne\p$ of $F$.
\item[(4)] If $\Psip(L)$ is a lattice, then
$\p^{2}\nn(L)\subset\nn(\Psip(L))\subset\p\cdot\nn(L)$.
\end{itemize}
The analogues of \textup{(1)} and \textup{(4)} hold for $\os_\p$-lattices.
\end{lem}

\begin{proof}
(1) We have $\nn(\p L)=\p^{2}\nn(L)\subset\p\nn(L)$, so $\p L\subset\Psip(L)$. The inclusion $\Psip(L)\subset L$ is strict because $\nn(L)$ is
generated by the values $Q(v)$, $v\in L$.

(2) Comparing $\ord_\q$ for all primes $\q$ of $F$, we see that
$\nn(L)\cap\p\nn(L)\os_\p=\p\nn(L)$. Since $Q(v)\in\nn(L)$ for every $v\in L$ and
$\nn(L)\os_\p=\nn(L_\p)$, it follows that a vector $v\in L$ satisfies $Q(v)\in\p\nn(L)$ if and
only if $Q(v)\in\p\nn(L_\p)$, i.e., if and only if $v\in\Psip(L_\p)$.

(3) Put $K:=L\cap\Psip(L_\p)$, so that $\Psip(L)=K$ by (2). As $\Psip(L_\p)$ is an
$\os_\p$-lattice by assumption, $K$ is an $\os$-module and $\p L\subset K\subset L$ by (1).
Hence $K$ is a lattice on $V$ with $K_\q=L_\q$ for all primes $\q\ne\p$ and $K_\p\subset\Psip(L_\p)$.
Conversely, let
$x\in\Psip(L_\p)$ and choose $v\in L$ with $v-x\in\p L_\p$, which is possible because $L$ is
dense in $L_\p$. Then
\[
Q(v)=Q(x)+2B(x,v-x)+Q(v-x)\in\p\nn(L_\p),
\]
since $Q(x)\in\p\nn(L_\p)$, $2B(x,v-x)\in2\p\sca(L_\p)\subset\p\nn(L_\p)$ by \eqref{eq:nsn},
and $Q(v-x)\in\p^{2}\nn(L_\p)$. Hence $v\in K$ and $x\in v+\p L_\p\subset K_\p$, because
$\p L\subset K$. This proves $K_\p=\Psip(L_\p)$.

(4) The right-hand inclusion is Definition \ref{def:psi}. The left-hand one follows from $\p
L\subset\Psip(L)$, since $\nn(\p L)=\p^{2}\nn(L)$.
\end{proof}

\begin{lem}\label{lem:E}
Let $L$ be an $\os$-lattice and let $\p$ be a prime of $F$ such that $\Psip(L)$ is a lattice. Then
$E(\Psip(L))\subset E(L)$. In particular, if $L$ is regular then so is $\Psip(L)$.
\end{lem}

\begin{proof}
Write $K:=\Psip(L)$. Since $\p L\subset K\subset L$ by Lemma \refpart{lem:psibasic}{(1)}, we have
$FK=FL$, $K_\q=L_\q$ for all primes $\q\ne\p$, and $K_\p\subset L_\p$, hence $Q(\gen K)\subset Q(\gen L)$. Let $\alpha\in E(K)$ and assume $\alpha=Q(v)$ for some $v\in
L$. As $Q(v)=\alpha\in Q(K_\p)\subset\nn(K_\p)\subset\p\nn(L_\p)$ by
Lemma \refpart{lem:psibasic}{(4)}, we get $v\in L\cap\Psip(L_\p)=K$ by Lemma \refpart{lem:psibasic}{(2)}, and
so $\alpha\in Q(K)$, a contradiction. Therefore $\alpha\notin Q(L)$ and $\alpha\in E(L)$.
\end{proof}

\begin{prop}\label{prop:index}
Let $\p$ be a prime of $F$ and let $K$ be an $\os_\p$-lattice of rank $n$ with $\nn(K)\supsetneq
2\sca(K)$. Then $\p$ is dyadic, $\Psip(K)$ is an $\os_\p$-lattice, and $[K:\Psip(K)]=\p$.
Consequently, if $L$ is an $\os$-lattice of rank $n$ with $\nn(L_\p)\supsetneq2\sca(L_\p)$, then
$\Psip(L)$ is an $\os$-lattice with
\[
[L:\Psip(L)]=\p,\qquad \vv(\Psip(L))=\p^{2}\vv(L),\qquad
\dd(\Psip(L))\supset\p^{-n+2}\dd(L).
\]
\end{prop}

\begin{proof}
From $2\sca(K)\subsetneq\nn(K)\subset\sca(K)$ we get $2\sca(K)\subsetneq
\sca(K)$, so $2\notin\uu_\p$ and $\p$ is dyadic, and \cite[Lemma 4.3(iv)]{Be} (where
$\Psip(K)$ is denoted by $K'$) shows that $\Psip(K)$ is an $\os_\p$-lattice with
$K/\Psip(K)\cong\os_\p/\p$, i.e., $[K:\Psip(K)]=\p$.

For the global statement, $\Psip(L)$ is a lattice by Lemma \refpart{lem:psibasic}{(3)}. As
$\Psip(L)_\q=L_\q$ for all primes $\q\ne\p$, we get $[L:\Psip(L)]=[L_\p:\Psip(L_\p)]=\p$ and
$\vv(\Psip(L))=\p^{2}\vv(L)$. Finally $\nn(\Psip(L))\subset\p\nn(L)$ gives
$\nn(\Psip(L))^{-n}\supset\p^{-n}\nn(L)^{-n}$, so
\[
\dd(\Psip(L))=2^{n}\nn(\Psip(L))^{-n}\vv(\Psip(L))\supset
\p^{-n+2}2^{n}\nn(L)^{-n}\vv(L)=\p^{-n+2}\dd(L).\qedhere
\]
\end{proof}

When $\nn(L_\p)=2\sca(L_\p)$, the set $\Psip(L)$ need not be a lattice, as in the example above,
but we will determine it when $L_\p$ has a suitable orthogonal splitting.

\begin{lem}\label{lem:split}
Let $\p$ be a prime of $F$ and let $K=K'\perp K''$ be an orthogonal splitting of an
$\os_\p$-lattice $K$ with $\nn(K')=\nn(K)$ and
$\nn(K'')\subset\p\cdot\nn(K)$. Then $\Psip(K)=\Psip(K')\perp K''$.
\end{lem}

\begin{proof}
For $x\in K'$ and $y\in K''$ we have $Q(x+y)=Q(x)+Q(y)$ with $Q(y)\in\nn(K'')\subset\p\nn(K)$.
Hence $Q(x+y)\in\p\nn(K)$ if and only if $Q(x)\in\p\nn(K)=\p\nn(K')$.
\end{proof}

\begin{lem}\label{lem:psimodular}
Let $\p$ be a prime of $F$ and let $M$ be a nonzero $\os_\p$-lattice such that
$Q(v)\os_\p=\nn(M)$ for every primitive $v\in M$. Then $\Psip(M)=\p M$.
\end{lem}

\begin{proof}
Let $v\in M$ be nonzero and write $v=\pi_\p^{c}u$ with $c\ge0$ and $u\in M$ primitive. Then
$\ord_\p Q(v)=2c+\ord_\p\nn(M)$, so $Q(v)\in\p\nn(M)$ if and only if $c\ge1$, i.e., if and
only if $v\in\p M$.
\end{proof}

\begin{prop}\label{prop:shape}
Let $L$ be an $\os$-lattice and let $\p$ be a prime of $F$. Let $L_\p=M\perp N$ be an
orthogonal splitting in which $\nn(M)=\nn(L_\p)\supsetneq\nn(N)$ and
$Q(v)\os_\p=\nn(M)$ for every primitive $v\in M$. Then $\Psip(L)$ is an $\os$-lattice with
$\Psip(L)_\p=\p M\perp N$.
\end{prop}

\begin{proof}
Since $\nn(N)\subsetneq\nn(M)=\nn(L_\p)$ we have $\nn(N)\subset\p\nn(L_\p)$, so
Lemmas \ref{lem:split} and \ref{lem:psimodular} give
$\Psip(L_\p)=\Psip(M)\perp N=\p M\perp N$. In particular $\Psip(L_\p)$ is an $\os_\p$-lattice,
so the assertion follows from Lemma \refpart{lem:psibasic}{(3)}.
\end{proof}

For dyadic $\p$, Lemma \ref{lem:split} is \cite[Lemma 4.3(i)]{Be}, and \cite[Lemma 4.10]{Be}
gives $\Psip(L_\p)=\p M\perp N$ when $M\cong\Ani^{\pi_\p^{i}}$, a case which we will use in
Section \ref{sec:local}.

\section{Reduction at a fixed prime}\label{sec:local}

In this section we will apply $\Psip$ repeatedly at a fixed prime $\p$. To show that this
process terminates, we will use the integer $\ord_\p\dd(K)$, where $K$ is the localization at $\p$
of the current lattice. Recall that $\dd(K)=2^{n}\nn(K)^{-n}\vv(K)$ for an $\os_\p$-lattice $K$ of
rank $n$. By Lemma \ref{lem:d}, $\ord_\p\dd(K)\ge0$.

\begin{prop}\label{prop:red}
Let $\p$ be a prime of $F$ and let $L$ be a regular $\os$-lattice of rank $n\ge4$. Then there is an
integer $t\ge0$ such that $\Psip^{t}(L)$ is a regular $\os$-sublattice of $L$ on $FL$ with
\[
Q\big(\Psip^{t}(L)_\p\big)=2\sca\big(\Psip^{t}(L)_\p\big).
\]
\end{prop}

\begin{proof}
For $i\ge0$ we put $L^{(i)}:=\Psip^{i}(L)$ whenever this is a lattice. It is then a lattice on
$FL$, since $\p^{i}L\subset L^{(i)}$ by Lemma \refpart{lem:psibasic}{(1)}. It suffices to show that,
as long as $Q(L^{(i)}_\p)\ne2\sca(L^{(i)}_\p)$, the lattice $L^{(i+1)}$ is again a regular
$\os$-sublattice of $L$, and that the nonnegative integers $\ord_\p\dd(L^{(i)}_\p)$ do not
increase with $i$ and decrease strictly within every two consecutive steps, because then there
is $t\ge0$ such that $\Psip^{t}(L)=L^{(t)}$ is a regular $\os$-sublattice of $L$ on $FL$ with
$Q(L^{(t)}_\p)=2\sca(L^{(t)}_\p)$.

We will first prove two claims about an $\os_\p$-lattice $K$ with $\nn(K)=2\sca(K)$ and
$Q(K)\ne2\sca(K)$. Let $K=M\perp N$, where $M$ is a leading component of $K$, and put
$\p^{a}:=\sca(M)=\sca(K)$.

\smallskip
{\bf Claim 1.} {\it The lattice $M$ is anisotropic, $\rank M\le2$, $Q(v)\os_\p=\nn(M)$ for every
primitive $v\in M$, and $\nn(N)\subsetneq\nn(M)=\nn(K)$.}

\smallskip
{\it Proof.} If $\p$ is non-dyadic, then $\nn=\sca$ for all $\os_\p$-lattices, so
$\nn(N)=\sca(N)\subsetneq\sca(M)=\nn(M)=\nn(K)$. If $M$ was isotropic, then
Proposition \refpart{prop:nondy}{(2)} would give
$\p^{a}=Q(M)\subset Q(K)\subset\nn(K)=\p^{a}=2\sca(K)$, a contradiction. Hence $M$ is
anisotropic, and Proposition \refpart{prop:nondy}{(3)} gives $\rank M\le2$ and
$Q(v)\os_\p=\p^{a}=\nn(M)$ for every primitive $v\in M$.

If $\p$ is dyadic, then $\nn(M)\subset\nn(K)=2\sca(K)=2\sca(M)\subset\nn(M)$, so
$\nn(M)=\nn(K)=2\sca(M)$.
By Corollary \refpart{cor:dycoro}{(1),(2)} and since $Q(K)\ne2\sca(K)$, the lattice $M$ is anisotropic
and $\nn(N)\ne\nn(M)$, so that $\nn(N)\subsetneq\nn(M)$. Proposition \refpart{prop:dymodular}{(2)}
then gives $\rank M=2$ and $Q(v)\os_\p=\nn(M)$ for every primitive $v\in M$. \qed

\smallskip
{\bf Claim 2.} {\it If $\rank M=\rank N=2$ and $\nn(N)=\p\nn(M)$, then
$\nn(\p M\perp N)\supsetneq2\sca(\p M\perp N)$, or $\p M\perp N$ has a leading component of
rank $1$.}

\smallskip
{\it Proof.} Put $K':=\p M\perp N$. If $\p$ is non-dyadic, then
$\sca(N)=\nn(N)=\p\nn(M)=\p^{a+1}$. The lattice $N$ is not modular, since otherwise
Corollary \ref{cor:nondy4} would give $Q(K)=\p^{a}=2\sca(K)$. Thus
$N=N_1\perp N_2$ with $\rank N_1=\rank N_2=1$ and $\sca(N_2)\subsetneq\sca(N_1)=\p^{a+1}$. As
$\sca(\p M)=\p^{a+2}$, the lattice $K'=N_1\perp\p M\perp N_2$ has a Jordan splitting with
leading component $N_1$.

If $\p$ is dyadic, then $\nn(N)=\p\nn(M)=2\p^{a+1}$ by Claim 1. Let us show that
$\nn(N)\ne2\sca(N)$. Otherwise a leading component $M_1$ of $N$ would satisfy
$2\p^{a+1}=2\sca(M_1)\subset\nn(M_1)\subset\nn(N)=2\p^{a+1}$, so that
$Q(N)\supset Q(M_1)\supset2\pi_\p^{a+1}\uu_\p\os_\p^{2}$ by
Proposition \refpart{prop:dymodular}{(1)}. Corollary \refpart{cor:dycoro}{(3)} would then give
$Q(K)=2\sca(K)$, a contradiction. Hence $2\sca(N)\subsetneq\nn(N)=2\p^{a+1}$ by \eqref{eq:nsn},
so $\sca(N)\subset\p^{a+2}=\sca(\p M)$. Consequently $\sca(K')=\p^{a+2}$ and
$\nn(K')=\nn(N)=2\p^{a+1}\supsetneq2\sca(K')$. \qed

\medskip
Now let $i\ge0$, assume that $L^{(i)}$ is a regular $\os$-sublattice of $L$, put
$K:=L^{(i)}_\p$, and assume that $Q(K)\ne2\sca(K)$. First, if $\nn(K)\supsetneq2\sca(K)$, then
Proposition \ref{prop:index} and Lemma \ref{lem:E} show that $L^{(i+1)}$ is a regular
$\os$-sublattice of $L$ with $\dd(L^{(i+1)})\supset\p^{-n+2}\dd(L^{(i)})$, i.e.,
\begin{equation}\label{eq:case1}
\ord_\p\dd\big(L^{(i+1)}_\p\big)\le\ord_\p\dd(K)-(n-2).
\end{equation}

Thus, by \eqref{eq:nsn}, it remains to consider the case $\nn(K)=2\sca(K)$. Let $K=M\perp N$,
where $M$ is a leading component of $K$, and put $r:=\rank M$. By Claim 1,
Proposition \ref{prop:shape}, and Lemma \ref{lem:E}, $L^{(i+1)}$ is a regular $\os$-sublattice
of $L$ with
\[
L^{(i+1)}_\p=\p M\perp N.
\]
As $\nn(N)\subset\p\nn(M)$, we have $\nn(\p M\perp N)=\p^{2}\nn(M)+\nn(N)=\p^{j}\nn(K)$, where
$j:=1$ if $\nn(N)=\p\nn(M)$ and $j:=2$ otherwise. Moreover $\vv(\p M\perp N)=\p^{2r}\vv(K)$
by \eqref{eq:perp} and \eqref{eq:scaling}. Therefore
\begin{equation}\label{eq:case2}
\ord_\p\dd\big(L^{(i+1)}_\p\big)=\ord_\p\dd(K)+2r-nj\le\ord_\p\dd(K),
\end{equation}
because $r\le2$, $j\ge1$, and $n\ge4$.

Equality holds in \eqref{eq:case2} only if $n=4$, $r=2$, and $j=1$, i.e., only if
$\rank M=\rank N=2$ and $\nn(N)=\p\nn(M)$. In that case, assume that
$Q(L^{(i+1)}_\p)\ne2\sca(L^{(i+1)}_\p)$. If $\nn(L^{(i+1)}_\p)\supsetneq2\sca(L^{(i+1)}_\p)$,
then \eqref{eq:case1} with $i+1$ in place of $i$ gives
$\ord_\p\dd(L^{(i+2)}_\p)<\ord_\p\dd(L^{(i+1)}_\p)$. Otherwise $L^{(i+1)}_\p$ has a leading
component of rank $1$ by Claim 2, and \eqref{eq:case2} with $i+1$ in place of $i$ and $r=1$
gives the same inequality.
\end{proof}

\section{Global results}\label{sec:main}

In this section we will prove Theorem \ref{thm:main}, rescale the resulting lattice $L'$, and
construct regular lattices of larger ranks.

\begin{proof}[Proof of Theorem $\ref{thm:main}$]
We will apply Proposition \ref{prop:red} successively over the set $S$ of primes
$\p$ of $F$ with $Q(L_\p)\ne2\sca(L_\p)$. This set is finite. Indeed, all but finitely many
primes $\p$ are non-dyadic and satisfy
$\vv(L_\p)=\os_\p$. For such $\p$ we have $\sca(L_\p)=\os_\p$ by Lemma \ref{lem:d}, so $L_\p$
is unimodular. Since $n\ge3$, the lattice $L_\p$ is isotropic by
Proposition \refpart{prop:nondy}{(1)}, and therefore $Q(L_\p)=\os_\p=2\sca(L_\p)$ by
Proposition \refpart{prop:nondy}{(2)}.

Write $S=\{\p_1,\dots,\p_r\}$ and put $L^{(0)}:=L$. For $u=1,\dots,r$, we apply
Proposition \ref{prop:red} to the lattice $L^{(u-1)}$ at the prime $\p_u$ and let
$L^{(u)}:=\Psi_{\p_u}^{t_u}(L^{(u-1)})$ be the resulting regular $\os$-sublattice of
$L^{(u-1)}$. Thus $FL^{(u)}=FL$ and $Q(L^{(u)}_{\p_u})=2\sca(L^{(u)}_{\p_u})$. By
Lemma \refpart{lem:psibasic}{(1)} we have $L^{(u)}_\q=L^{(u-1)}_\q$ for every prime $\q\ne\p_u$, so the
equality $Q(L^{(u)}_{\p_w})=2\sca(L^{(u)}_{\p_w})$ remains valid for all $w\le u$, and
$L^{(r)}_\p=L_\p$ for every prime $\p\notin S$. Hence $L':=L^{(r)}$ is a regular
$\os$-sublattice of $L$ with $FL'=FL$ and
\[
Q(L'_\p)=2\sca(L'_\p)\qquad\text{for every prime }\p\text{ of }F.
\]

Now let $\alpha\in\os$ be nonzero. Since $L'$ is positive definite, we have
$\sigma_i(\alpha)\in Q_{\sigma_i}((FL')_{\sigma_i})$ for all $i$ if and only if
$\alpha\in\os^{+}$.
Moreover $\alpha\in Q(L'_\p)=2\sca(L'_\p)$ for all primes $\p$ if and only if
$\alpha\in2\sca(L')$, because $\sca(L')\os_\p=\sca(L'_\p)$ for all $\p$. Therefore
$Q(\gen L')\setminus\{0\}=2\sca(L')\cap\os^{+}$. As $L'$ is also regular,
we get
\[
Q(L'\setminus\{0\})=Q(L')\setminus\{0\}=Q(\gen L')\setminus\{0\}=2\sca(L')\cap\os^{+}.\qedhere
\]
\end{proof}

Replacing the lattice $L'$ from Theorem \ref{thm:main} by $(L')^{\alpha}$ for suitable totally
positive $\alpha\in F$ gives the following.

\begin{prop}\label{prop:narrowclass}
Let $F$ be a totally real number field and let $L$ be a positive definite regular $\os$-lattice
of rank $n\ge4$. Then there is a narrow ideal class $C$ of $F$ such that for every integral
ideal $\mf b\in C$ there is a positive definite regular $\os$-lattice $M$ of rank $n$ with
$\sca(M)=\mf b$ and
\[
Q(M\setminus\{0\})=2\mf b\cap\os^{+}.
\]
\end{prop}

\begin{proof}
Let $L'$ be as in Theorem \ref{thm:main}, put $\mf a:=\sca(L')$, which is a nonzero integral
ideal, and let $C$ be the narrow ideal class of $\mf a$. Let $\mf b\in C$ be integral and write
$\mf b=\alpha\mf a$ with $\alpha\in F^{+}$. Put $M:=(L')^{\alpha}$, which has rank $n$ and is
positive definite because $\alpha$ is totally positive. By \eqref{eq:scaling} we have
$\sca(M)=\alpha\mf a=\mf b\subset\os$, so $M$ is integral and hence regular by
Lemma \ref{lem:scaleE}. Finally, multiplication by $\alpha$ maps $2\mf a\cap\os^{+}$
bijectively onto $2\mf b\cap\os^{+}$, since $\alpha(2\mf a)=2\mf b$ and
$\alpha F^{+}=F^{+}$, so that
\[
Q(M\setminus\{0\})=\alpha Q(L'\setminus\{0\})=\alpha\big(2\mf a\cap\os^{+}\big)
=2\mf b\cap\os^{+}. \qedhere
\]
\end{proof}

\begin{coro}\label{cor:narrow}
Let $F$ be a totally real number field of narrow class number $1$ and let $L$ be a positive
definite regular $\os$-lattice of rank $n\ge4$. Then there is a positive definite regular
$\os$-lattice $L''$ of rank $n$ with $\sca(L'')=\os$ and $Q(L''\setminus\{0\})=2\os^{+}$.
\end{coro}

\begin{proof}
Since the narrow class number of $F$ is $1$, the narrow ideal class $C$ from
Proposition \ref{prop:narrowclass} contains $\os$. Proposition \ref{prop:narrowclass} with
$\mf b=\os$ therefore gives the lattice $L''$, because $2\os\cap\os^{+}=2\os^{+}$.
\end{proof}

For $\eta=(\eta_1,\dots,\eta_d)\in\{\pm1\}^{d}$ put
$\os^{(\eta)}:=\{\alpha\in\os\setminus\{0\}\mid \sgn(\sigma_i(\alpha))=\eta_i
\text{ for all }1\le i\le d\}$, so that $\os^{(1,\dots,1)}=\os^{+}$.

\begin{coro}\label{cor:class}
Let $F$ be a totally real number field of class number $1$ and let $L$ be a positive definite
regular $\os$-lattice of rank $n\ge4$. Then there are $\eta\in\{\pm1\}^{d}$ and a
regular $\os$-lattice $L''$ of rank $n$ with $\sca(L'')=\os$ and
\[
Q(L''\setminus\{0\})=2\os^{(\eta)}.
\]
\end{coro}

\begin{proof}
Let $L'$ be as in Theorem \ref{thm:main}. Since the class number of $F$ is $1$, there is
$\alpha\in\os$ with $\sca(L')=\alpha\os$. Put
$L'':=(L')^{\alpha^{-1}}$ and
$\eta:=(\sgn\sigma_1(\alpha),\dots,\sgn\sigma_d(\alpha))$. Then $\sca(L'')=\os$ and, by
\eqref{eq:scaling},
\[
Q(L''\setminus\{0\})=\alpha^{-1}\big(2\alpha\os\cap\os^{+}\big)
=2\{\beta\in\os\setminus\{0\}\mid \alpha\beta\in\os^{+}\}=2\os^{(\eta)},
\]
because $\sgn\sigma_i(\alpha\beta)=\eta_i\sgn\sigma_i(\beta)$ for every $i$. By
Lemma \ref{lem:scaleE}, $L''$ is regular.
\end{proof}

Finally, we will show that a totally real number field that admits a positive definite regular
lattice of rank $n\ge4$ also admits such a lattice of every rank $m\ge n$. 

\begin{prop}\label{prop:rankup}
Let $F$ be a totally real number field and let $L$ be a positive definite regular $\os$-lattice of
rank $n\ge4$. Then there is $a\in\Z^+$ such that $L\perp\qf{a}$ is regular.
\end{prop}

\begin{proof}
We will choose $a$ such that $Q(\gen(L\perp\qf{a}))=Q(\gen L)$.
Let $S$ be the finite set of primes $\p$ of $F$ that are dyadic or satisfy
$\vv(L_\p)\ne\os_\p$. As in the proof of Theorem \ref{thm:main}, we have $Q(L_\p)=\os_\p$ for
every prime $\p\notin S$.

Let $\p\in S$. Since the quadratic space $F_\p L_\p$ has dimension $n\ge4$, it represents every
element of $F_\p^{\times}$ by \cite[63:21]{O}. For each of the finitely many classes $c$ of
$F_\p^{\times}$ modulo $F_\p^{\times2}$ we choose $x_c\in L_\p$ with $Q(x_c)\in c$. Put
$k_\p:=\max_c\ord_\p Q(x_c)$, so that $k_\p\ge0$. If $\beta\in\p^{k_\p}$ is nonzero and lies in
the class $c$, then $\beta=\mu^{2}Q(x_c)$ for some $\mu\in F_\p^{\times}$ with
$2\ord_\p\mu=\ord_\p\beta-\ord_\p Q(x_c)\ge0$, hence $\beta=Q(\mu x_c)\in Q(L_\p)$. Thus
\begin{equation}\label{eq:allsmall}
\p^{k_\p}\subset Q(L_\p)\qquad(\p\in S).
\end{equation}

Choose a positive rational integer $a$ with $a\in4\p^{k_\p}$ for every $\p\in S$ and put
$M:=L\perp\qf{a}$. Then $M$ is a positive definite $\os$-lattice containing $L$, so
$Q(\gen L)\subset Q(\gen M)$.

Conversely, let $\alpha\in Q(\gen M)$. The conditions at the
real embeddings $\sigma_i$ are the same for $L$ and $M$, since both lattices are positive
definite. For every prime $\p\notin S$ we have $\alpha\in\os_\p=Q(L_\p)$. Let $\p\in S$ and write
$\alpha=Q(x)+at^{2}$ with $x\in L_\p$ and $t\in\os_\p$. If $\alpha\in\p^{k_\p}$, then
$\alpha\in Q(L_\p)$ by \eqref{eq:allsmall}. Otherwise $\alpha\ne0$ and $\ord_\p\alpha<k_\p$,
hence
\[
\ord_\p\big(at^{2}/\alpha\big) \ge \ord_\p a-\ord_\p\alpha \ge \ord_\p(4)+1,
\]
so $1-at^{2}/\alpha\in1+4\p$ is the square of a unit $u$ by the local square theorem
\cite[63:1]{O}. Thus $\alpha=Q(u^{-1}x)\in Q(L_\p)$ in this case as well. We proved that
$\alpha\in Q(\gen L)$, so $Q(\gen M)=Q(\gen L)$. As $L$ is regular,
\[
Q(\gen M)=Q(\gen L)=Q(L)\subset Q(M)\subset Q(\gen M),
\]
so $M$ is regular.
\end{proof}

\section{Application to real quadratic fields}\label{sec:density}

In this section we will deduce Theorems \ref{thm:GRH} and \ref{thm:finite} from
Theorem \ref{thm:main}.

The \textit{absolute norm} of a nonzero integral ideal $\mf b$ is $N\mf b:=\#(\os/\mf b)$. As in the
introduction, $h^{+}(D)$ denotes the narrow class number of $\mathbb Q(\sqrt D)$ for a
squarefree integer $D>1$.

For an
integer $m\ge4$ we put
\[
\mathcal R_m:=\{D>1 \text{ squarefree}\mid \mathbb Q(\sqrt D) \text{ admits a positive definite
regular lattice of rank } m\}.
\]
By Proposition \ref{prop:rankup}, $D\in\mathcal R_m$ if and only if $\mathbb Q(\sqrt D)$ admits a
positive definite regular lattice of some rank $r$ with $4\le r\le m$.

For sets $\mathcal E\subset\mathcal D$, the {\it density of $\mathcal E$ in
$\mathcal D$} is
\[
\lim_{X\to\infty}\frac{\#\{D\in\mathcal E\mid D\le X\}}{\#\{D\in\mathcal D\mid D\le X\}}
\]
whenever this limit exists. The {\it lower density of $\mathcal E$ in $\mathcal D$} is the
corresponding limit inferior. For functions $f$ and $g$ of a real variable $X$ we write $f\ll
g$, or equivalently $g\gg f$, if $|f(X)|\le Cg(X)$ for some constant $C>0$ and all $X$ under
consideration. A subscript indicates the parameters on which $C$ is allowed to depend.

From Proposition \ref{prop:narrowclass} we obtain lattices of the kind counted in \cite{KYZ}.

\begin{coro}\label{cor:muuniv}
Let $k\ge1$ be a real number and let $F$ be a totally real number field such that every narrow
ideal class of $F$ contains an integral ideal of norm at most $k$. If $F$ admits a
positive definite regular $\os$-lattice of rank $n\ge4$, then there is
$\nu\in\Z^+, \nu\le k,$ and a positive definite $\os$-lattice of rank $n$ that represents every totally
positive multiple of $2\nu$.
\end{coro}

\begin{proof}
Let $C$ be the narrow ideal class provided by Proposition \ref{prop:narrowclass}, let $\mf b\in
C$ be an integral ideal with $\nu:=N\mf b\le k$, and let $M$ be the corresponding lattice. Then
$\nu\in\mf b$ because $\nu$ annihilates $\os/\mf b$, so $2\nu\os^{+}\subset2\mf
b\cap\os^{+}=Q(M\setminus\{0\})$.
\end{proof}

\begin{thm}\label{thm:density}
Let $m\ge4$ be an integer, put $\delta:=\frac1{48m}$, and let $\mathcal D$ be a set of
squarefree integers $D>1$ such that for every $D\in\mathcal D$, every narrow ideal
class of $\mathbb Q(\sqrt D)$ contains an integral ideal of norm at most $D^{\delta}$. Then
\[
\#\{D\in\mathcal D\cap\mathcal R_m\mid D\le X\} \ll_m X^{11/12}(\log X)^{3/2}.
\]
\end{thm}

\begin{proof}
Write $T(X):=\#\{D\in\mathcal D\cap\mathcal R_m\mid D\le X\}$. For positive
integers $R$ and $\mu$ let $B(R,\mu)$ be the quantity defined after \cite[Theorem 3.1]{KYZ}, in which we replace the determinant $\det G$ by $1$. Let
$\mathcal N(R,\mu,X)$ be the number of squarefree integers $1<D\le X$ for which there is a
positive
definite lattice of rank $R$ over $\mathbb Q(\sqrt D)$ representing every totally positive
multiple of $\mu$. By \cite[Theorem 4.2]{KYZ},
\[
\mathcal N(R,\mu,X)<100B(R,\mu)^{3/2}X^{7/8}(\log X)^{3/2}
\qquad\text{whenever } X\ge B(R,\mu)^{12}(\log X)^{4}.
\]
The explicit formula for $B(m,\mu)$ gives a constant $c=c(m)$ with
\begin{equation}\label{eq:Bbound}
B(m,\mu)\le c\mu^{m}\qquad (\mu\ge1).
\end{equation}

Let $D\in\mathcal D\cap\mathcal R_m$ with $D\le X$. Corollary \ref{cor:muuniv}, applied with
$k=D^{\delta}$, gives 
$\nu\in\Z^+, \nu\le X^{\delta},$ such that $\mathbb Q(\sqrt D)$ admits a positive definite lattice of rank
$m$ representing every totally positive multiple of $2\nu$. Consequently
\[
T(X) \le \sum_{1\le\nu\le X^{\delta}}\mathcal N(m,2\nu,X).
\]
By \eqref{eq:Bbound} we have $B(m,2\nu)\le c(2X^{\delta})^{m}$ for all the $\nu$ occurring
here. As $12m\delta=\frac14$, it follows that
$B(m,2\nu)^{12}(\log X)^{4}\le c^{12}2^{12m}X^{1/4}(\log X)^{4}\le X$ once $X$ exceeds a bound
depending only on $m$. For such $X$ the displayed bound of \cite[Theorem 4.2]{KYZ} applies to
each of the at most $X^{\delta}$ summands and yields
\[
T(X) \le 100c^{3/2}2^{3m/2}X^{7/8+(3m/2+1)\delta}(\log X)^{3/2}.
\]
Since $(\frac{3m}{2}+1)\delta\le2m\delta=\frac1{24}$ and $\frac78+\frac1{24}=\frac{11}{12}$, this
is the asserted bound.

The choice of $\delta$ only makes the exponents simple. The same argument works for every
$\delta<\frac1{12m+8}$ and gives the exponent $\frac78+(\frac{3m}{2}+1)\delta<1$ in place of
$\frac{11}{12}$.
\end{proof}

\begin{proof}[Proof of Theorem $\ref{thm:GRH}$]
 Let $\delta:=\frac1{48m}$ and $0<\eta<\delta/2$. The main task is to show that there is $D_0$
such that for every squarefree $D>D_0$ with $h^{+}(D)\le D^{\eta}$, every narrow ideal class of
$\mathbb Q(\sqrt D)$ contains an integral ideal of norm at most $D^{\delta}$.

Let $D>1$ be squarefree, let $F:=\mathbb Q(\sqrt D)$, let $\Delta_F\le4D$ be the discriminant of
$F$, and let $E$ be the narrow Hilbert class field of $F$, i.e., the maximal abelian extension of
$F$ that is unramified at all finite primes. By class field theory, the Artin map induces an
isomorphism from the narrow ideal class group of $F$ onto $\mathrm{Gal}(E/F)$ that maps the class
of a prime ideal of $F$ to its Artin symbol. In particular, $E$ has degree $2h^{+}(D)$ over
$\mathbb Q$. As $E/F$ is unramified at all finite primes, the absolute discriminant of $E$ is
$\Delta_E=\Delta_F^{h^{+}(D)}$. Under the generalized Riemann hypothesis, \cite[Theorem 3.1(1)]{BS}, applied to $E/F$, therefore
shows that every narrow ideal class of $F$ contains a prime ideal of norm at most
\begin{equation}\label{eq:BS}
	(1+o(1))\big(\log\Delta_E+2[E:\mathbb Q]\big)^{2}=(1+o(1))\big(h^{+}(D)(\log\Delta_F+4)\big)^{2}.
\end{equation} If
$h^{+}(D)\le D^{\eta}$, then this bound is 
less than $D^{\delta}$ for all sufficiently large $D$. This proves the existence of $D_0$.

Let $\mathcal D$ be the set of squarefree integers $D>D_0$ with $h^{+}(D)\le D^{\eta}$. As we
have just shown, Theorem \ref{thm:density} applies to $\mathcal D$, and so $\mathcal D\cap\mathcal R_m$
has density $0$ in the set of all squarefree integers. By the second assumption of
Theorem \ref{thm:GRH}, the complement of $\mathcal D$ has density $0$ as well, and
hence so does $\mathcal R_m$.
\end{proof}

\begin{rmk}\label{rmk:hooley}
For squarefree $D>1$, let $h(D)$ be the class number of $\mathbb Q(\sqrt D)$. Hooley's
conjecture \cite[Conjecture 7]{Hoo} immediately implies
$\sum_{\text{squarefree }D\le X}h(D)\le X^{1+o(1)}$. Hence for every $\eta>0$, counting the
$D\le X^{1/2}$ trivially and using $h(D)\ge h^{+}(D)/2>X^{\eta/2}/2$ for the other $D$, we get
\[
\#\{D\le X\mid D \text{ squarefree},\ h^{+}(D)>D^{\eta}\} \le X^{1/2}+2X^{-\eta/2}
\sum_{\text{squarefree }D\le X}h(D) = o(X).
\]
Thus the second assumption of Theorem \ref{thm:GRH} follows from \cite[Conjecture 7]{Hoo}.
\end{rmk}

The Cohen--Lenstra heuristics \cite{CL} predict that $\mathbb Q(\sqrt p)$ has class number $1$
for approximately $75.4\%$ of the primes $p\equiv1\pmod 4$, and the same is expected for all
primes \cite{tRW}, where numerical evidence for $p\equiv1\pmod 4$ is also given. In particular, the heuristics predict that the assumption of the following corollary holds.

\begin{coro}\label{cor:CL}
Let $m\ge4$ be an integer. Assume that the set of primes $p\not\equiv7\pmod 8$ for which
$\mathbb Q(\sqrt p)$ has class number $1$ has positive lower density in the set of all primes.
Then
\[
\#\{p\le X\mid p \text{ prime},\ p\notin\mathcal R_m\} \gg \frac{X}{\log X}.
\]
\end{coro}

\begin{proof}
Let $\mathcal D$ be the set of primes $p\ge2^{48m}$ with $p\not\equiv7\pmod 8$ for which
$\mathbb Q(\sqrt p)$ has class number $1$. We claim that for every $p\in\mathcal D$, every
narrow ideal class of $\mathbb Q(\sqrt p)$ contains an integral ideal of norm at most $2$.

If $p\equiv1\pmod4$, then the fundamental unit of $\mathbb Q(\sqrt p)$ has norm $-1$, so the
narrow class number of $\mathbb Q(\sqrt p)$ equals its class number $1$, and the claim holds.

Let $p\equiv3\pmod8$. Then $\os=\Z[\sqrt p]$ and no unit of $\os$ has norm $-1$, since $-1$ is
not a square modulo $p$. As all ideals of $\os$ are principal, it follows that the narrow ideal
classes are the class of $\os$ and the class of the ideals $\alpha\os$ with $\alpha$ of negative
norm. The
prime ideal of $\os$ above $2$ has norm $2$ and is generated by some $\alpha=x+y\sqrt p$ with
$x,y\in\Z$ and $x^{2}-py^{2}=\pm2$. Here $x^{2}-py^{2}=2$ is impossible, since $2$ is not a
square modulo $p$, so this prime ideal lies in the second class. This proves the claim.

As $2\le p^{1/(48m)}$ for all $p\in\mathcal D$, Theorem \ref{thm:density} applies to
$\mathcal D$ and gives $\#\{p\in\mathcal D\cap\mathcal R_m\mid p\le X\}=o(X/\log X)$. The set
$\mathcal D$ contains all but finitely many of the primes in our assumption, so
$\#\{p\in\mathcal D\mid p\le X\}\gg X/\log X$, and we get $\gg X/\log X$ primes $p\le X$ in
$\mathcal D\setminus\mathcal R_m$.
\end{proof}

\begin{proof}[Proof of Theorem $\ref{thm:finite}$]
Let $F=\mathbb Q(\sqrt D)$ have narrow class number $1$ and admit a positive definite regular
lattice of rank $r$ with $4\le r\le7$. By Corollary \ref{cor:narrow} there is a positive
definite regular $\os$-lattice $L''$ of rank $r$ with
$Q(L''\setminus\{0\})=2\os^{+}$, so the
asserted bounds are those of \cite[Theorem 1.2]{KK+} for lattices of rank at most $r$ that
represent every element of $2\os^{+}$.
\end{proof}

Without the assumption on the narrow class number, a similar argument gives the following under
the generalized Riemann hypothesis.

\begin{coro}\label{cor:GRHclass}
Assume the generalized Riemann hypothesis. If $D>1$ is a squarefree
integer such that $\mathbb Q(\sqrt D)$ admits a positive definite regular lattice of rank $r$
for some $r$ with $4\le r\le7$, then
\[
h^{+}(D)\gg\frac{D^{1/32}}{(\log D)^{3/2}}.
\]
\end{coro}

\begin{proof}
We will first choose constants that do not depend on $D$. Let $\Delta_K$ denote the
discriminant of a real quadratic field $K$. By \cite[Theorem 1.2]{KK+} under the generalized
Riemann hypothesis, there are $M_0\ge2$ and $C>0$ such that $\Delta_K\le C\mu^{16}(\log\mu)^{16}$
for every integer $\mu\ge M_0$ and every real quadratic field $K$ that admits a positive definite
lattice of rank at most $7$ representing every totally positive multiple of $\mu$. By the
unconditional bounds of \cite[Theorem 1.2]{KK+}, there is $D_1$ such that $\Delta_K<D_1$
whenever such a lattice exists for some integer $\mu$ with $2\le\mu<M_0$. Finally, we choose
$D_2\ge D_1$ such that for all $D>D_2$, the factor $1+o(1)$ in \eqref{eq:BS} is at most $2$ and
$2\big(D^{1/32}(\log4D+4)\big)^{2}\le D$.

Now let $D$ be as in the statement, put $F:=\mathbb Q(\sqrt D)$, so that $D\le\Delta_F\le4D$, and assume that
$h^{+}(D)\le D^{1/32}$, since otherwise there is nothing to show. For $D\le D_2$ the assertion
holds after enlarging the implied constant, so let $D>D_2$. By \eqref{eq:BS} and the choice of
$D_2$, every narrow ideal class of $F$ contains an integral ideal of norm at most
\[
k:=2\big(h^{+}(D)(\log4D+4)\big)^{2},
\]
and $k\le D$. Corollary \ref{cor:muuniv} therefore gives $\nu\in\Z^+,\nu\le k,$ and a positive
definite lattice of rank $r$ over $F$ that represents every totally positive multiple of $2\nu$.
Since $\Delta_F>D_1$, we have $2\nu\ge M_0$ by the definition of $D_1$, and so
\[
D\le\Delta_F\le C(2\nu)^{16}(\log2\nu)^{16}.
\]
As $\nu\le k\le D$, this yields $D\ll h^{+}(D)^{32}(\log D)^{48}$.
\end{proof}

We conclude with a few remarks on our results. The restriction to \textbf{ranks at least $4$} in
Theorem \ref{thm:main} cannot be removed in general, since by the argument of \cite[Lemma 3]{EK}, its conclusion is false for all ternary
lattices over every totally real field of odd degree. Over fields of even degree, our proof of
Theorem \ref{thm:main} nevertheless uses $n\ge4$ in \eqref{eq:case2}, to ensure that
$\ord_\p\dd(\Psip^{i}(L)_\p)$ does not increase with $i$.

\medskip

In the proofs above, we used \textbf{narrow ideal classes} because the lattices counted in \cite{KYZ} are positive
definite. However, the results of \cite{KYZ} can be extended to lattices $L$ of mixed signature (i.e., $Q(L)\subset\os^{(1,-1)}\cup\{0\}$), which, together with Corollary \ref{cor:class}, would then allow us to drop the condition
$p\not\equiv7\pmod8$ in Corollary \ref{cor:CL}.

Moreover, as $h^{+}(D)$ differs from the class number $h(D)$ by a factor of $1$ or $2$, the
asymptotic results (Theorem \ref{thm:GRH} and Corollary \ref{cor:GRHclass}) can equally well be
formulated with $h(D)$ in place of $h^{+}(D)$.

\medskip

Finally, recall that all our lattices are \textbf{integral} (or classical), i.e., $B(L,L)\subset\os$. Our results also
cover {\it integer-valued} (or {\it non-classical}) lattices $L$, for which only
$Q(L)\subset\os$ is required. For such a lattice $L$ we have
$2B(v,w)=Q(v+w)-Q(v)-Q(w)\in\os$ for all $v,w\in L$, so $L^{2}$ is an integral lattice of the same rank. If $L$ is positive definite and regular, then
so is $L^{2}$, by the proof of Lemma \ref{lem:scaleE}. Thus Theorems
\ref{thm:main}--\ref{thm:finite} and Corollaries \ref{cor:CL} and \ref{cor:GRHclass} hold for
integer-valued lattices as well.

\end{document}